\documentclass[hidelinks,onefignum,onetabnum]{siamart250211} 

\usepackage[margin=1in,letterpaper]{geometry}
\usepackage{lipsum}
\usepackage{amsfonts}
\usepackage{graphicx}
\usepackage{epstopdf}
\usepackage{algorithmic}
\ifpdf
  \DeclareGraphicsExtensions{.eps,.pdf,.png,.jpg}
\else
  \DeclareGraphicsExtensions{.eps}
\fi

\newsiamremark{remark}{Remark}
\newsiamremark{hypothesis}{Hypothesis}
\crefname{hypothesis}{Hypothesis}{Hypotheses}
\newsiamthm{claim}{Claim}
\newsiamremark{fact}{Fact}
\crefname{fact}{Fact}{Facts}
\newcommand{\ds}{\displaystyle}
\headers{ROM and Anderson type acceleration for ill-posed problems}{K. Ito and T. Xue}

\title{Reduced order method based Anderson-type acceleration method for least square problems and large scale ill-posed problems \thanks{Submitted to the editors \today.
}} 

\author{Kazufumi Ito\thanks{Department of Mathematics, Graduate Program of Operations Research, Center for Research in Scientific Computation, North Carolina State University, Corresponding author.
  (\email{kito@ncsu.edu}, \url{https://kito.wordpress.ncsu.edu/}).}
\and Tiancheng Xue \thanks{Grauate Program of Operations Research, North Carolina State University
  (\email{txue2@ncsu.edu}).}
} 

\usepackage{amsopn}

\newcommand{\R}{\mathbf{R}}
\newcommand{\Z}{\mathbf{Z}}

\newcommand{\half}{\frac{1}{2}}

\ifpdf
\hypersetup{
  pdftitle={Anderson},
  pdfauthor={Kazufumi Ito， Tiancheng Xue}
}
\fi

\begin{document}

\maketitle

\begin{abstract}
In this paper, we propose an acceleration framework for a class of iterative methods 
using the Reduced Order Method (ROM). Assuming that the underlying iterative scheme 
generates a rich basis for the solution space, we construct 
the next iterate by minimizing the equation error over the linear manifold spanned by this basis.
The resulting optimal linear combination yields a more accurate approximation of the solution 
and significantly enhances convergence.  In essence, the method can be seen as a history-based
acceleration technique, akin to a delayed or memory-enhanced iterative scheme. This approach 
effectively remedies semi-ill-posed problems, enabling convergence where standard methods may fail, 
and also acts as a stabilizing and regularizing mechanism for the original iteration.

\end{abstract}



\section{Introduction}
Anderson acceleration is a widely recognized technique known for its effectiveness in enhancing the convergence of fixed-point problems. By leveraging information from previous iterations, it combines weighted past iterations to generate a new one. Originally introduced to accelerate the iteration process of nonlinear integral equations, the method has since been extended to general fixed-point problems \cite{WN}. Anderson acceleration has been extensively studied and applied, yielding promising results in practical computations,
e.g., \cite{A}\cite{WN} \cite{SS} \cite{TXHSX} and references therein.

\textbf{\underline{Contributions}}: In this paper, we extend Anderson acceleration and apply a general iterative method to the nonlinear least square problem and ill-posed problems. In particular, we formulate and analyze it as the reduced order method \cite{IR} \cite{WI} \cite{IK2008} \cite{An} to achieve a monotone convergence. Our method incorporates a given iterative method as subspace generation process. As a consequence, the use of Anderson acceleration remedies iterative algorithms where relative residual may fail to converge under tolerance level. In addition, our method works on a broader class of problems. For example, we simplify conjugate gradient method with variable step size iterative method and make it work on an indefinite system of linear equations, as we shall see from our presentation.

\textbf{\underline{Notations}}: Throughout the paper, let $X$ be a Banach space with inner product structure $(\cdot, \cdot)$ and norm $|\cdot|_2$, and $A$ be a linear operator over $X$.

A central theme of the paper is to consider designing an algorithm that can accelerate the convergence of iterative solutions over a linear manifold spanned by solutions generated from a given iterative scheme. For example, we consider variable step size iterative method 
\begin{equation}\label{VSSI}
x_{n+1}=x_n+\alpha_np_n,
\end{equation}
where $\{x_n\}$ is a sequence of solution update, and $p_n$ is a search direction, say $p_n=F(x_n)$, for solution $x_n$ in $X$ with step size $\alpha_n$ on
a large scale optimization problem of minimizing 
$$
|F(x)|^2+\gamma\,\Theta(x),
$$
where $F:X\to X$ is semi smooth map and $\Theta$ defines the regularization of unknown
$x \in X$  by the Anderson-type method. This encompass many applications. 
For example, given  transport map $T_\theta$,
one can formulate the least squares problem on the optimal transport problem \cite{CIY} \cite{DF} \cite{Itma}:
\begin{equation} \label{opt}
|T_\theta(g_0)-g_1|^2+\Theta(\theta) \mbox{  over $\theta$}, 
\end{equation}
where $\Theta$ is the transport cost. 

\noindent\underline{\textbf{Example}: Network transport} 
Let $\{x_i\}$ be the nodes in graph $G$ with nodes $S$ and edge $E$, set 
$$
L_{e,i}=-1; L_{e,j}=1,\; \mbox{   if $x_i \in S \to x_j \in S$ is connected by edge $e \in E$}, 
$$
and let a symmetric positive definite matrix $D$ be the matrix of edge weight. Now, with net source vector $s$, one can formulate a quadratic transport problem as
$$ \min \quad \half \theta^{\top} D^{-1} \theta \ \text{subject to} \ L\theta= s. $$ 
Put
$
A=L^{\top}DL.
$
The Lagrange multiplier can imply that putting $T_{\theta}(g_0) = A \theta$ and $g_1 =s$ defines the transport problem \eqref{opt}, regardless of $g_0$. 

\noindent The plan of our presentation is as follows:
\begin{itemize}
    \item We set up the foundation of acceleration method. We examine various aspects of improvement 
and extension of
the proposed iterative methods in Sections \ref{asp}, e.g., variable step size,  subspace iterate, Quasi Newton, (conjugate) gradient for nonlinear optimization. At the same time, we introduce extrapolation as method for acceleration in Section \ref{extp} and incorporate it with schemes examined in Section \ref{asp}. 
\item Afterwards, we develop and analyze the acceleration method
based on the reduced order method (ROM) \cite{IR}. We introduce ROM and Anderson acceleration in Section \ref{exts}. We interpret Anderson as generating evolutional fixed point \cite{Ka} \cite{CP}. Our key observation is 
that iterative methods naturally generate a solution subspace structure that can serve as the ROM basis to 
accelerate the underlying iteration. By combining a well-constructed solution basis 
with the ROM framework, we obtain the desired acceleration effect.
We further examine the connection between ROM and Anderson acceleration. 
In the nonlinear setting, Anderson acceleration produces an approximation that is closely related to the ROM solution, 
and the ROM optimality criterion can be used as a merit function for guiding the Anderson method. 
\item Finally, we detail on applications and the corresponding analysis for
specific examples. Specifically in Section \ref{StkT} we examine the saddle point problem on Stokes type system. Afterwards, we refine a specific iterative iterative scheme to work on in Section \ref{VssFp}.
We then present convergence analysis and numerical findings, on matrix Riccati system in Section \ref{Experiment3} and
Stokes type system in Section \ref{Exper4}. 
\end{itemize}
\section{Various interactive methods}\label{asp} 
In this section we discuss interactive methods for generating
solution basis. Specifically, we discuss variable step size iterative method \eqref{VSSI}:
$$
x_{n+1}=x_n+\alpha_np_n.
$$
\subsection{Fixed point problem: $F(x)=0$}\label{FPP}
We formulate merit function minimization as an extension of least square problem for fixed point of another map. Suppose we want to find a fixed point for $G: X\to X$. We can consider minimizing $|x-G(x)|^2$ over $x\in X$. Let $F(x)= x-G(x)$, then we have minimizing $|F(x)|^2$. With $f^\prime(x) = F(x)$, we have defined an equivalent problem for finding a fixed point for $G$: $\min\;\; f(x)$. 
With $J\approx f''(x_n)$ by BFGS \cite{BFGS}, $J^{-1}$ computed by the Sherman-Morrison-Woodbury \cite{SM} \cite{W} formula, and $f^\prime$ upgraded by Broyden method \cite{B2}, we have: 
$$ x_{n+1}-x_n=\alpha\, J^{-1}f'(x_n).$$
\subsection{Quasi Newton}
We first discuss Newton method. We assume the (damped) Newton updates $u_n$. For $1\le n\le m$,
$$
u_{n+1} =u_n -J_n^*(J_nJ^*_n)^{-1}F(u_n),
$$
where $J_n$ is an approximated Jacobian as $F^\prime(u_n)$.


In a quasi-Newton iteration, 
$$
x_{n+1} =x_n +\alpha_n B_n^{-1}F(x_n),
$$
where $B_k$ is some approximation of the Jacobian matrix $F^\prime(x_k)$, e.g., Broyden update.
For example the Barzilai-Borwein method \cite{BB}: $B\sim \beta I$ with
$$
\beta_k=\frac{(\Delta x,\Delta x)}{(\Delta x,\Delta r)}\mbox{  or  } 
\beta_k=\frac{(\Delta r,\Delta x)}{(\Delta r,\Delta r)},
$$
where
$$
\Delta x=x_{n}-x_{n-1} , \Delta r=F_{n}-F_{n-1}.
$$
In the context of feedback control \cite{BI}, we can consider the feedback form:
$$
x_{n+1} =x_n-\alpha_nG_nF(x_n),
$$
where $G_n$ is a feedback gain operator.
\subsection{Conjugate gradient}
In this section we consider the conjugate gradient method and ROM acceleration. For a symmetric positive matrix $A$
we consider
$$
\min \quad \frac{1}{2}(x,Ax)-(b,x).
$$
The necessary optimally condition is $Ax-b=0$. We start with  $r_0 = b-Ax_0, p_0 = r_0$. Then, for $1\le k \le m$, the pre-conditioned conjugate gradient method is given by
$$\begin{array}{l}
\alpha_k= \frac{(r_k,z_k)}{(Ap_k,p_k)}, \;\; x_{k+1}=x_k+\alpha_k\,p_k,
\\ \\
r_{k+1}=r_k-\alpha_k\, Ap_k,\;\; z_{k+1}=P^{-1}r_{k+1}, 
\\ \\
\beta_k=\frac{(z_{k+1},r_{k+1}),}{(z_k,r_k)},\;\; p_{k+1}=z_{k+1}-\beta_k p_k. 
\end{array} $$
Here, $P$ is the pre-conditioner and $z_k$ is the pre-conditioned residual. If the CG method starts with $x_0 = 0$,  then $x_{k}$ minimizes
$$
x_k =\mbox{argmin} \{(x-y,A(x-y)): y \in \mbox{span}(b,Ab,\cdots,A^{k-1}b)\}.
$$

\subsection{Nonlinear conjugate gradient}
For the nonlinear convex minimization
$$
\min\;\; f(x) \mbox{  subject to } x \in {\cal C},
$$
where ${\cal C}$ is a closed convex set. For $1\le n \le m$,
The conjugate gradient method is
$$
x_{n+1}-x_n=\alpha_n\,p_n,\;\;p_n=f^{\prime}(x_n).
$$
The step size $\alpha_n$ minimizes
$$
f(x_n+\alpha p_n),\mbox{  subject to  }x_n+\alpha p_n\in {\cal C},
$$
where
$$
p_{n+1}=r_{n+1}-\beta_n\, p_n.
$$
With $\Delta x_n = -\nabla_x f(x_n)$, for the conjugacy
we use the Polak–Ribiere formula:
\begin{equation}\label{Polak-Ribiere F}
\beta^{\text{PR}}_n = \frac{\Delta x_n^{\top} (\Delta x_n -\Delta x_{n-1})}{\Delta x_{n-1}^{\top} \Delta x_{n-1}},    
\end{equation}
instead of the Fletcher–Reeves formula: 
\begin{equation}\label{Fletcher-Reeves F}
    \beta^{\text{FR}}_{n} = \frac{\Delta x_n^{\top} \Delta x_n}{\Delta x_{n-1}^{\top} \Delta x_{n-1}}.
\end{equation} 
\begin{remark}
But, one can improve the convergence of above methods, say CG, with ROM acceleration. Specifically, one can modify CG method so that it can be applied without assuming $A$ being positive definite and symmetric with a modified step size $\alpha_k$, say Cauchy step with ROM acceleration.    
\end{remark}

\section{Solution extrapolation}\label{extp}
We introduce DIIS as the method for updating solution through Anderson acceleration. DIIS \cite{P} updates the current solution as a linear combination from previous iterations. 
\begin{theorem}
    Given basis vectors $(x_1, \cdots, x_m) \in X \times \cdots \times X$, with $e_i = F(x_i)$ for $i \in \{ 1,\cdots,m \}$, consider the problem:
\begin{equation}\label{DIIS}
        \min \quad |\sum_{i=1}^m c_i e_i|^2 \ \text{subject to}\ \sum_i c_i = 1.
\end{equation}
With $B_{ij} = ( e_i , e_j )$, the optimal solution $c_1, \cdots, c_n$ is obtained by solving the linear equation
    \begin{equation}\label{DIISSys}
    \begin{bmatrix}
        B_{11} & \cdots & B_{1m} & 1\\
        \vdots & \ddots & \vdots & \vdots \\
        B_{m1} & \cdots & B_{mm} & 1\\
        1 & \cdots & 1 & 0
    \end{bmatrix}
    \begin{bmatrix}
        c_1\\
        \vdots\\
        c_m\\
        \lambda
    \end{bmatrix}
    =  \begin{bmatrix}
        0 \\
        \vdots\\
        0 \\
        1
    \end{bmatrix}.
\end{equation}
\end{theorem}
\begin{proof}
With $B_{ij} = ( e_i , e_j )$, one formulate the Lagrangian \cite{IK} associated to \eqref{DIIS} as:
    $$ L(c,\lambda) = c^{\top} B c - 2\lambda (\sum_i c_i -1).$$
    Taking partial derivatives of $L$ with respect to the coefficients and the multiplier and equating to 0 leads to: 
    \begin{equation}\label{DIISSysO}
    \begin{bmatrix}
        B_{11} & \cdots & B_{1m} & -1\\
        \vdots & \ddots & \vdots & \vdots \\
        B_{m1} & \cdots & B_{mm} & -1\\
        1 & \cdots & 1 & 0
    \end{bmatrix}
    \begin{bmatrix}
        c_1\\
        \vdots\\
        c_m\\
        \lambda
    \end{bmatrix}
    =  \begin{bmatrix}
        0 \\
        \vdots\\
        0 \\
        1
    \end{bmatrix}.
    \end{equation}
Equating minus sign on $\lambda$ and treating $- \lambda$ as $\lambda$ results in \eqref{DIISSys}.
\end{proof}
Therefore, as a consequence, to obtain Anderson update, we do linear combination of previous iterate with coefficients of \eqref{DIISSys} by
\begin{equation}\label{ROMUpdate}
x_{m+1} = \sum_{i=1}^{m} c_i x_i.     
\end{equation}

\begin{remark}
Given $\epsilon > 0$, one can include a regularization sub-matrix $\epsilon I_m$ into the \eqref{DIISSys} so that one solves \begin{equation}\label{DIISSysReg}
    \begin{bmatrix}
        B_{11} + \epsilon & \cdots & B_{1m} & 1\\
        \vdots & \ddots & \vdots & \vdots \\
        B_{m1} & \cdots & B_{mm}+\epsilon & 1\\
        1 & \cdots & 1 & 0
    \end{bmatrix}
    \begin{bmatrix}
        c_1\\
        \vdots\\
        c_m\\
        \lambda
    \end{bmatrix}
    =  \begin{bmatrix}
        0 \\
        \vdots\\
        0 \\
        1
    \end{bmatrix},
    \end{equation}
    to update the variable coefficients for the desired system with weight coefficients $c_1, \cdots, c_m$.
\end{remark}
\section{ROM, Anderson acceleration, and
Evolutional Fixed point} \label{exts}
The main goal of this section is to first define ROM and Anderson type acceleration, then describe the relations between them, and finally develop variants of them, such as Nested ROM, Sampled ROM, and ROM extrapolation, which will be essential for our numerical experiments in Section \ref{Experiment3} and Section \ref{Exper4}. ROM construct the next iterate by minimizing the equation error over the linear manifold spanned by the basis generated by the underlying iterative scheme. Anderson minimizes the total sum of residual. In addition, we interpret Anderson as an evolutional fixed point method. We use ROM in bringing convenience in theoretical analysis, while Anderson acceleration is better for practice for giving extrapolation solution. 
\subsection{ROM in general} In this section, we provide ROM formulations for merit function and residual function.

\begin{definition}
Given $f: \R^n \to \R$ as a merit function, ROM step we minimize
$$
\min \quad f(\sum_k\alpha_k x_k) \mbox{  subject to  }\sum_k\alpha_k x_k\ \in {\cal C}
$$
and let
$$
x_{m+1}=\sum_{k \le m} \alpha^*_k\,x_k.
$$    
\end{definition} 
\begin{definition}
Given $F: X \to X$ in residual form, ROM step gives us the convergence based on the criterion:  $\alpha^*_k$ minimizes
\begin{equation} \label{ROM}
\min\quad |F(\sum_{k\le m}\alpha_k\,x_k)|^2, 
\end{equation}
and update by
$$
x_{m+1}=\sum_{k\le m} \alpha^*_k x_k.
$$    
\end{definition} 
\noindent \textbf{\underline{Example}}: When $F(x) = b-Ax$ and $r_k = F(x_k)$, we have:
$$|F(\sum_{k\le m} \alpha_kx_k)| = |b - A(\sum_{k\le m} \alpha_kx_k) |^2. 
$$ 
\subsection{Anderson type acceleration}
In this section we develop  the acceleration  of Anderson type for the standard fixed point  iterate
$$
x_{i+1}=g(x_{i}).
$$ 
That is, let  $r_k=x_k-g(x_k)=F(x_k)$, solve for $\alpha$ to 
\begin{equation}\label{anderson}
\min\quad| \sum_{k\le m}\alpha_k\,r_k|^2 \mbox{  subject to  } \sum_{k\le m} \alpha_k=1,
\end{equation}
and then the sequential update $x_{m+1}$ is defined by
$$
x_{m+1}=\sum_{k\le m} \alpha_k^*\,x_k,
$$
where $\alpha^*$ is an optimizer of \eqref{anderson}.  It is modified and different from the original form of the original Anderson update \cite{A}. The step regularizes and stabilizes the fixed point update. It works without the contraction of fixed point iterate. That is, it defines an evolutionary
fixed point with an improved contraction. For the evolutional fixed point problem, we let $\{u_i\}$ be fixed point solution to one parameter family of 
damped fixed points and consider:
$$
u_i-\theta_i\, g(u_i)=0,\;\; |\theta_i|<1.
$$
and $u_i$ is the corresponding fixed point for the damped fixed point problem.



\subsection{Comparison of Anderson and ROM}\label{AndROM}
In this section we compare ROM and Anderson.
\begin{proposition}
If $F: X \to X$ defined by $F: x \to b-Ax$ with $r_k =F(x_k)$ is linear, consider constrained ROM, i.e. \eqref{ROM} with the constraint $\sum_{k\le m}\alpha_k=1$, ROM is equivalent to minimizing the total residual
$$
|\sum_{k\le m}\alpha_k\ r_k|^2 \;\;  \mbox{  subject to } \sum _{k\le m}\alpha_k=1. 
$$    
\end{proposition}
\begin{proof}
Notice that $\sum_{k\le m} \alpha_k = 1$ implies $\sum_{k\le m} \alpha_k \ b = b$. Further, linearity of $A$ suggests
    \begin{align*}
       |F(\sum_{k\le m} \alpha_kx_k)| 
       &= |\sum_{k\le m} \alpha_k \ (b - A x_k) |^2 := |\sum_{k\le m}\alpha_k\ r_k|^2 .
    \end{align*}
Therefore, Anderson provides the same solution as constrained ROM, given a linear residual function. 
\end{proof}
\begin{theorem}
    Let $F: X\to X$, consider 2 problems: the first one is \eqref{anderson}, i.e.,  
$$
\min\quad|\sum_{k\le m}\alpha_k\,r_k|^2 \mbox{  subject to  } \sum_{k \le m} \alpha_k=1.
$$ 
and \eqref{ROM} with the constraint $\sum_{k\le m}\alpha_k=1$. Then, for some function $H: X\to X$ such that $\ds \lim_{x_k \to \bar{x}} H(x_k) = 0$ for every $\bar{x}$ close to $x_k$, we have:
$$ \sum_k \alpha_k r_k \sim F(\sum_k \alpha_k x_k) + H(x_k),$$
\end{theorem}
\begin{proof}
Note that given $\bar{x}$ close to $x_k$, there exists $H: X\to X$ such that $\ds \lim_{x_k \to \bar{x}} H(x_k) = 0,$ and $$ F(x_k) \sim F(\bar{x}) + F'(\bar{x}) (x_k - \bar{x}) + H(x_k).$$
After adding and subtracting $F(\bar{x})$ and $\sum_k \alpha_k \bar{x} = \bar{x}$, we do linear approximation as:
\begin{align*}
    \sum_k \alpha_k r_k &= \sum_k \alpha_k (F(x_k) - F(\bar{x})) + F(\bar{x})\\
    &\sim \sum_k \alpha_k F'(\bar{x}) (x_k - \bar{x}) + F(\bar{x}) + H(x_k)\\
    &\sim F'(\bar{x}) (\sum_k \alpha_k x_k - \sum_k \alpha_k \bar{x}) +  F(\bar{x}) + H(x_k)\\
    &\sim F(\bar{x}) + F'(\bar{x}) (\sum_k \alpha_k x_k - \sum_k \alpha_k \bar{x})  + H(x_k)\\
    &\sim F(\sum_k \alpha_k x_k) + H(x_k).
    \label{approx3}
\end{align*} 
\end{proof}
That is, the Anderson method gives an approximate solution that corresponds to the first order approximation of constrained ROM:
$$
\min_{\alpha} \quad |F(\sum_k \alpha_k \,x_k)|^2 \mbox{  subject to } \sum_k \alpha_k=1,
$$
and with $J \sim F'(\bar{x})$, $x = \sum \alpha_k x_k$, and $\bar{x} = \sum \bar{\alpha}_k x_k$, we have:
\begin{equation} \label{GN}
\min_{\alpha} \quad |\sum_k  J\,(\alpha_k-\bar{\alpha}_k)\,x_k+F(\bar{x})|^2.
\end{equation}

\subsection{Nested ROM}\label{Nested}
We perform the nested ROM, i.e., we progressively generate sample solution sequence $\{z_i\}$ with $1\le i \le m$ and $m \in \Z_{\ge 2}$, by the Anderson updates. 
ROM with $m$ step can be nested as: let $y_1 = y $ be $2$-step ROM solution be defined by
$$
y=\sum^2_{k=1} \alpha^*_k x_k,
$$
with $\alpha^*_k \in \R$ as optimized weight for iterate $x_k$. Then, based on the previous solutions, we generate an additional $r$ based on the initial condition $y$ and use 3-step ROM to get $y_2$, so on. Then, after $m$ times, we progressively obtain a sampled solution sequence $\{y_i\}$. We obtain a refined approximation by nested ROM update
$$
z=\sum^m_{k=1} \beta^*_k y_k,
$$
where $\beta^*=(\beta^*_1, \cdots, \beta^*_m)\in \R^m$ minimizes
$$
|F(\sum^m_{k=1}\beta^*_k y_k)|^2.
$$
Also, we do an additional refinement based on $\{z_i\}_{1 \le i \le m}$, as a nested ROM. It notes that this improves accuracy and the convergence very much and that they produce a monotone descent approximation sequences for minimizing  $|F(x)|^2.$


It follows the reduced order element method \cite{IR}
with the iterates $\{x_k\}$ as basis, i.e., the linear manifold spanned by $\{x_k\}$ captures the solution well by minimizing
\begin{equation}
    \phi(\sum_{k\le m}\alpha_k x_k) = |F(\sum_{k\le m}\alpha_k x_k)|^2 + \beta\,\Psi(\sum_{k\le m}\alpha_k x_k),
\end{equation}
where $\beta\ge 0$ and $\Psi$ is a regularization function.

As an application, we discuss Nested Anderson with an application on improving the fixed point algorithm. Suppose that we want to find a fixed point for a given contractive map $\phi: \R^n \to \R^n$. Define $F(x) = x - \phi(x)$. Let the Anderson method define a fixed point map $\phi$ by
\begin{equation}\label{NestROM}
    \phi(\theta)=\theta^*,
\end{equation}
where $\theta$ is the initial and $\theta^*$ is the corresponding the $m$-step Anderson iterate. In general, the Anderson map \eqref{NestROM} has a much better contraction than $\phi$. That is, we apply the Anderson method  for the fixed point problem for $\phi$ as a nested iterate.

\subsection{Sampled ROM}
Sampled ROM is a similar method to Nested ROM except to restart after each $m$-step ROM. We repeat $m$-step ROM for sometime with update and we collect the update together as a new basis for our new ROM. Mathematically speaking, Sampled ROM step can be described as: let $y_1 = y $ be $m$-step ROM solution be defined by
$$
y_1=\sum^m_{k=1} \alpha^*_k x_k,
$$
with $\alpha^*_k \in \R$ as optimized weight for iterate $x_k$. Then, we restart ROM from $y_1$ and repeat the $m$-step ROM method with initial $y$ to get $y_2$, so on. Then, we progressively obtain a sampled solution sequence $\{y_i\}$. We obtain a refined approximation by nested ROM update
$$
z=\sum^m_{k=1} \beta^*_k y_k,
$$
where $\beta^*=(\beta^*_1, \cdots, \beta^*_m)\in \R^m$ minimizes
$$
|F(\sum^m_{k=1}\beta_k y_k)|^2.
$$
\subsection{ROM extrapolation}\label{StkT}
Consider the least square problem of
$$
\min\quad |F(x)|^2.
$$
We introduce the ROM update
$$
\min \quad |F(\sum_{k\le m}\,\alpha_k x_k)|^2,
$$
where the basis 
$(x_k,r_k)$ is generated by  a well-posed  linear system  $\bar{A}x-b=0$
$$
x_{n+1}=x_n+\alpha_n\, p_n, \;\; r_n=b-\bar{A}x_n.
$$
That is, we generate basis functions
by the alternative nominal equation  $\bar{A}x-b=0$ for 
the original  (ill-posed) problem
of minimizing $|F(x)|^2$.  For example, equation $F(x)=0$ is a perturbation of $\bar{A}x-b=0$. Let $F(x)$ be a linear equation of the form
$$
F(x)=(\bar{A} + \tilde{A})x - b. 
$$ 
For example, consider the saddle point problem
$$
\bar{A} = \left(\begin{array}{cc} -D & E^* \\ \\ E &O\end{array}\right), \tilde{A} =\left(\begin{array}{cc} \delta F & O^* \\ \\ O &O\end{array}\right), A =\left(\begin{array}{cc} -D+\delta F & E^* \\ \\ E &O\end{array}\right).
$$
For the case of the Navier Stokes we have
$$
\bar{A}=\mbox{the Stokes operator and  }
\delta F=\vec{u}\cdot\nabla \vec{u}.
$$
We formally introduce Saddle Point Problem in Section \ref{Saddle}.

\section{Variable step size fixed point}\label{VssFp}
We develop an acceleration method
for various iterative methods in numerical optimizations and
variational problems. We discuss a general class of
Interactive methods and  their acceleration of the convergence.
We use the reduced order method (ROM)
acceleration. It results  in a rapidly convergent algorithm. The proposed method can be applied to many applications.
\begin{definition} 
Let $F: X\to X$ be a Lipschitz map on a Banach space $X$. The Cauchy step $\alpha_n \in \R$ at $x_n \in X$ along direction $p_n \in X$ is defined as $$\alpha_n^* = \min |F(x_n + \alpha p_n)|^2 \ \text{over $\alpha\in \R$.}$$ 
\end{definition}
We now explain the variable fixed point iterate through the use of Cauchy Step. For minimizing $|F(x)|^2$ 
, we have the variable fixed point iterate 
$$
x_{k+1}=x_k+\alpha_k\,F(x_k),\;  p_k=F(x_k),
$$
for the fixed point problem $F(x)=x-g(x)=0$,
where $\alpha_n$ is a variable step size: since the step size $\alpha_n$  preconditioning of each iterate $\{x_n\}$, it is very essential to find an optimal value for $\alpha_n$ for the convergence and stability of the iterative method.
In the case of the fixed point we have $F (x) = x-g(x)$ and it defines the variable fixed point iterate. For example the variable fixed point for minimizing $ |Ax -b|^2$ we have
\begin{equation} \label{L}
x_{n+1} =x_n +\alpha_n\,(b-Ax_n) ,
\end{equation}
where $\alpha_n$ is a variable step size. 
\begin{theorem}\label{CTU}
The optimal step size for \eqref{L} is Cauchy step given by
\begin{equation}\label{ACUno}
\alpha_n = \frac{(r_n,Ap_n)}{(Ap_n,Ap_n)}.
\end{equation}    
\end{theorem}
\begin{proof}
Set $G(\alpha) = |b - A(x_n + \alpha_n d_n)|_2^2$, evaluate $\frac{dG}{d\alpha_n}$, and set $\frac{dG}{d\alpha_n} = 0$ to get $$(b - A(x_n + \alpha_n d_n))^{\top} (Ad_n) = 0.$$ With re-arrangement, one gets $$ \alpha_n = \frac{(b - Ax_n,Ad_n)}{(Ad_n, Ad_n)} =\frac{(r_n, Ad_n)}{(Ad_n, Ad_n)}.$$
In the context of variable fixed point, we have $d_n = r_n$, then we have: $$ \alpha_n = \frac{(r_n, Ar_n)}{(Ar_n, Ar_n)},$$  as desired.
\end{proof}
\begin{corollary}
Let $F: X\to X$ be nonlinear, consider finding the optimal step size $\alpha$ that minimizes
$$
|F(x+\alpha\,p)|^2, 
$$
we have:
$$
\alpha_n =-\frac{(r_n,Jp_n)}{(Jp_n,Jp_n)},
$$
where
$$
Jp=\frac{F(x+yp)-F(x)}{t},\;\; t>0.
$$
\end{corollary}
One can also consider the subspace update spanned by multi-search directions of $\{q_j\}$.
$$
x_{k+1} =x_k+\sum_j\beta_j\, q_j.
$$
For example, we consider the subspace spanned by 
$(r_k,F^\prime(x_k)r_k)$:
$$
x_{k+1} =x_k+\alpha_k r_k+\beta_k F^\prime(x_k)r_k.
$$
For each $j$, we minimize over step sizes $\beta_j$
$$
|F(x_k+\sum_j \beta_j\,q_j)|^2. 
$$
It is a subspace method via $span \{q_j\}$. As a consequence of theorem \ref{CTU}, we have the following:
\begin{corollary}
Consider the subspace spanned by  $(r_k,Ar_k)$:
\begin{equation} \label{sub}
x_{k+1} = x_k+ (\alpha_k r_k+\beta_k Ar_k),  
\end{equation}
where $(\alpha_n, \beta_n)$ is generated by Cauchy step:
\begin{equation}\label{CDos}
\left(\begin{array}{c}
     \alpha_n  \\
     \\
     \beta_n 
\end{array} \right)=
\left(\begin{array}{cc}
(Ar_n, Ar_n) & (Ar_n,Ap_n) \\ \\
(Ap_n, Ar_n) &(Ap_n, Ap_n)\end{array} \right)^{-1}
\left(\begin{array}{c}
    (r_n, Ar_n) \\
    \\
    (r_n, Ap_n)
\end{array}
\right).
\end{equation}    
\end{corollary}



Thus, we sequentially generate $\{x_1,\cdots, x_m\}$ by \eqref{sub} -
\eqref{CDos} and use it as a basis for the solution. That is,
we use the reduced order  method (ROM) \cite{IR}
$$
\min\;\;
|F(\sum_k \gamma_k x_k)|^2+\beta\,\psi(\sum_k \gamma_k x_k)
$$
and set the extrapolation  of $\{x_k\}$ by
\begin{equation}
x_{m+1}=\sum_{k \le m} \gamma_k ^*x_k.
\end{equation}
It also introduces the stabilization and regularization of iterates 
$\{x_k\}$
and we select $\psi$ to enhance such properties. The ROM step accelerates the convergence and stabilizes and regularizes the method. It is an optimal polynomial acceleration and  the Lagrange extrapolation of $(x_k,F(x_k))$. However, not at all time will the variable step size update work, we use the idea of ROM extrapolation, and propose an alternative Cauchy step. 

\noindent \underline{An alternative Cauchy Step} to minimize $|b-Ax|^2$ over $x\in X$: We use ROM extrapolation introduced in Section \ref{StkT}. Let $\bar{A}$ be an approximated matrix for $A$. Then, we set
\begin{equation}\label{ACUno}
    \alpha_n = \frac{(r_n, \bar{A}r_n)}{(\bar{A}r_n, \bar{A}r_n)}. 
\end{equation}
We update $x_k$ by 
$$
x_{k+1}=x_k+\alpha_k r_k, \quad r_k = b - Ax_k.
$$  
Similarly, as an alternative, for step size update, we have:  
\begin{equation}\label{ACDos}
\left(\begin{array}{c}
     \alpha_n  \\
     \\
     \beta_n 
\end{array} \right)=
\left(\begin{array}{cc}
(\bar{A}r_n, \bar{A}r_n) & (\bar{A}r_n,\bar{A}p_n) \\ \\
(\bar{A}p_n, \bar{A}r_n) &(\bar{A}p_n, \bar{A}p_n)\end{array} \right)^{-1}
\left(\begin{array}{c}
    (r_n, \bar{A}r_n) \\
    \\
    (r_n, \bar{A}p_n)
\end{array}
\right).
\end{equation}
We replace $A$ by an approximated matrix $\bar{A}$ in the matrix for inversion. For updating directions, we consider directions $r_n$, and $\bar{A}r_n$, with $r_n = b-Ax_n$:
\begin{equation} 
x_{k+1} = x_k+ (\alpha_k r_k+\beta_k Ar_k). 
\end{equation}
\section{Convergence} We bring the results from basis solution generation and ROM method together to establish convergence arguments.
Throughout this section, to make analysis easier, we let $x_{i,j}$ be the $m$-step update Anderson update, with $1 \le i\le m$ and $j \ge 0$. We re-establish conclusions from Section \ref{VssFp} with new notations first.

Accordingly, for the linear case $F(x) = b - Ax$, we update our notation for variable step update: for fixed $j$, and $1\le i < m$, we have the variable step update:
\begin{equation}\label{CVSU}
x_{i+1,j} = x_{i,j} + \alpha_{i,j}\ p_{i,j},
\end{equation}
we have the Cauchy step:
\begin{equation} \label{CUno}
\alpha_{i,j} =\frac{(r_{i,j},Ap_{i,j})}{(Ap_{i,j},Ap_{i,j})}.
\end{equation}
Similarly, with regularization function $|x|^2$, we have:  
$$
|Ax-b|^2+\beta\,|x|^2,\;\; \beta\ge 0,
$$
we have variable step update \eqref{CVSU} with Cauchy Step:
$$ \alpha_{i,j} = \frac{ r_{i,j}^{\top} Ap_{i,j} - \beta x_{i,j}^{\top} p_{i,j} }{|Ap_{i,j}|^2 + \beta |p_{i,j}|^2}, $$
that minimizes:
$$
|A(x_{i,j} +\alpha_{i,j} \ p_{i,j})-b|^2 +\beta\, |x_{i,j} +\alpha_{i,j} \ p_{i,j}| ^2.\;\; 
$$
\begin{theorem} Consider the least squares problem
$$
F(x) = |Ax-b|^2+\beta\,|x|^2,\;\; \beta\ge 0,
$$
For fixed $j$, the sequence $\{ x_{i,j} \}$ updated by variable step size update gives a monotonic decreasing sequence $\{|F(x_{i,j})|^2\}_{i=1}^m$.
\end{theorem}
\begin{proof} Fix $j$.
We first deal with the case $\beta = 0$, we have the following: $$
|F(x_{i+1,j})|^2=|F(x_{i,j})|^2-\frac{|(r_{i,j},Ap_{i,j})|^2}{(Ap_{i,j},Ap_{i,j})},
$$
Similarly, if $\beta \not = 0$, we have:  $$
|F(x_{i+1,j})|^2 = |F(x_{i,j})|^2 -\frac{( \beta x_{i,j}^{\top} p_{i,j} - r_{i,j}^{\top}Ap_{i,j})^2}{|Ap_{i,j}|^2 + \beta |p_{i,j}|^2}.
$$
Therefore, the sequence $\{|F(x_{i,j})|^2\}_{i=1}^m$ gives a monotonically decreasing sequence, with $\beta \ge 0$.
\end{proof}
We next establish results on $\{|F(x_{1,j})|^2\}_{j=0}^{\infty}$. Therefore, we first recast with notations in this section: 
let $g:X \to X$ with $F(x) =x -g(x)$, ROM is given by
\begin{equation} \label{ROM1}
\min\;\; |F(\sum _{k\le m } \alpha_{k,j} x_{k,j})|^2, 
\end{equation}
with update
$$
x_{1,j+1}=\sum_{k\le m} \alpha^*_{k,j} x_{k,j},
$$
where $\{\alpha^*_k\}$ is an optimizer. For convenience, put $r_{1,j} = F(x_{1,j})$.

\begin{lemma}
When $G$ is linear, we have
    \begin{equation}\label{linear}
r_{1,j+1}=\sum_{k\le m}\alpha^*_{k,j} r_{k,j}. 
\end{equation}    
\end{lemma} 
\begin{proof}
    Let $\alpha^*$ be optimal weight for the Anderson update. Note that 
    \begin{align*}\label{linear}
r_{1,j+1} &= x_{1,j+1}-g(x_{1,j+1}) \\
    &= \sum_{k\le m}\alpha^*_{k,j} x_{k,j} - \sum_{k\le m}\alpha^*_{k,j} g(x_{k,j})\\
    &= \sum_{k\le m}\alpha^*_{k,j} (x_{k,j} -  g(x_{k,j})) \\
    &= \sum_{k\le m}\alpha^*_{k,j} r_{k,j}.
\end{align*}
\end{proof}
\begin{theorem}\label{CL}
Given a linear $F$, the sequence $\{x_{1,j}\}_{j=0}^{\infty}$ generated by \eqref{ROM1} 
converges weakly to the optimal pair $(x^*,r^*)$, with $r^* = F(x^*)$. 
\end{theorem}
\begin{proof}
WLOG, we assume that $0 < \max\{\alpha_{i,j}\} < 1$ for all $1\le i \le m$ and $j$. So, there exists $\rho \in \R$ such that $0 < \rho < \max\{\alpha_{i,j}\} < 1$. As a result, we have:
    \begin{align*}
        |r_{1,j+1} | &\le \max_i\{\alpha_i\} \max_i |r_{i,j}|\\
        &\le \rho \max_i |r_{i,j}|\\
        & \le \rho^2 \max_i |r_{i,j-1}| \\
        & \cdots\\
        & \le \rho^{j+1} \max_i |r_{i,0}|.
    \end{align*}
    This suggests a contraction and implies convergence for $\{|r_{1,j}|\}$. 
\end{proof}

\begin{remark}
For the linear case, ROM is equivalent to the Krylov subspace method \cite{Kr} \cite{WN}.
Thus, its performance 
depends on the conditioning of the residual system $\{r_k\}$. If $\dim(X)$ is finite and the condition number of $A$ is not large, span $\{r_i\}$ does not increase then 
$(x_k,r_k) \to (x^*,r^*)$ in a finite step.  If 
$\{r_k\}$ is nearly singular, i.e., the ROM method slowly convergent if 
$\{r_k\}$ is badly conditioned. e use other methods, say extrapolation method. 
\end{remark}

\begin{remark} In general $r^*\neq 0$ and $x^*$ minimizes $|F(x)|^2$. When solving $Ax = b$, If $A\bar{x} = 0$ for some nonzero $\bar{x}$ and $b \notin R(A)$, then the pair $(x^*,r^*)$  is not unique.
Let $A = USV^*$ be a singular value decomposition $A$.  Assume that for 
$V = V_1 + V_2$ and
$$
|AV_2| \sim 0 \mbox{  and   }\min |AV_1- b| > 0 
$$
Then, it is numerically singular and divergent in the sense that the iterate very slowly
convergent to $(x^*,r^*),\; r^*\neq  0$. It happens if the condition number of $A$ is large. 
\end{remark}

\begin{remark}
For the nonlinear convex case, after making linearization, with same argument as linear case, one can conclude that there is a subsequence  
$\{(x_k, r_k)\} \to (x^*, r^*)$
weakly. If $F$ is weakly continuous, i.e. if $\{x_k\}$ converges to $x^*$,
then $\{F(x_k)\} \to F(x^*)=r^*$ weakly, then
$$
x^*-g(x^*)=r^*.
$$    
\end{remark}
For $F: X\to X$ being nonlinear, in order to establish the strong convergence,
we assume the coercivity
$$
(x-y-(g(x)-g(y)),x-y)\ge\omega |x-y|^2,
$$
then
$$
\omega\,|x_k-x^*| \le |r_k-r^*| \to 0  \mbox{ as  $k\to\infty$},
$$
and
$$
\omega\, |x_{i+1}-x_i| \le  |F(x_{i+1})-F(x_i)|.
$$
One can prove the following theorem by applying the same argument of Theorem \ref{CL} when $A$ is replaced by approximated Jacobian as $F'(x_n)$:
\begin{theorem}
Given a weakly continuous $F$, the sequence $\{x_{1,j}\}_{j=0}^{\infty}$ generated by \eqref{ROM1} 
converges weakly to the pair $(x^*,r^*)$. 
\end{theorem}
\begin{remark}
    Under the coercivity, if $|r^*|=0$ then $x_n \to x^*$, the unique fixed point. 
\end{remark}




\noindent \underline{\textbf{Example}: Operator equation}
We consider the operator equation of the form:
$Lx+ f(x) = 0$
where $L$ is boundedly invertible and $f$ is Lipschitz on $X$.
It Is formulated as a fixed point problem
$$
F(x)=x-L^{-1}f(x)=0
$$
and the variable fixed point
$$
x_{n+1}=x_n-\alpha_n F(x_n).
$$
Assume
$$
(Lx -f (x) -(Ly-f (y)), x- y) \ge  \omega\, |x -y|^2.
$$
Then, for some $c > 0$,
$$
\omega\,|x_i-x^*|\le c\,|r_i -r^*|.
$$
For example,
$$
L = -\Delta,\;\; f(u) = -\nabla\cdot (\phi^\prime(|\nabla u|^2)\nabla u),
$$
where $\phi$ is convex, with $X= H^1_0(\Omega)$.
Thus,
$$
(Lu+f(u)-(Lv+f(v)),u-v)\ge(L(u-v),u-v), 
$$
and thus,
$$
|x_n- x^*|_X \le  c |r_n-r^*|_{X^*},\; 0\le  c<1.
$$ 
\section{Matrix Riccati equation} \label{Experiment3} 
 In this section, we demonstrate the use of Anderson-type acceleration method for two specific examples and objectives. Also, we refer to extensive test examples in the literature, e.g.,\cite{A} \cite{WN}.
 
\subsection{Nonlinear system for matrix Riccati equation} 
 In this section, we demonstrate the use of Anderson map and continuation fmethod on a matrix Riccati equation. A matrix Ricccati equation behaves like Navier Stokes equation, i.e., quadratic nonlinearity. Specifically, we consider a matrix Riccati equation of the form 
    \begin{equation}\label{MRE}
         \dot{u} = -b\ u \ u + Q,
    \end{equation}
where $u$ is an unknown matrix of size $n\times n$, $Q$ is a symmetric positive matrix of $n \times n$, and $b >0$ is the internal variable for the map. It is a nonlinear system of quadratic equation and is parameter dependent with parameter $b$. Here, $b>0$ is like the Reynolds number in in the Navier Stokes system So, there are many issues associated with $b$. For example, large $b$ significantly magnifies nonlinearity, which brings challenges to find the matrix equation solution. Therefore, we are concerned about developing computational methods to solve a matrix Riccati equation. We solve \eqref{MRE} by performing the fixed point iterate by$F(u) = -b\ u \ u + Q$ and make sure $F$ is contractive, i.e., $b>0$  sufficiently small. Starting with small $b$, we assume that the map $g(u)=-b\ u\ u+Q$ is contractive, $u$ is square matrix in $R^{n\times n}$ and $u\ u$ denotes the matrix product. We see that large $b$ stops contraction after contraction range. We use Anderson map with the method of continuation to construct  solutions to large $b>0$.   The fixed point method does not converge for a critical $b$.
But, Anderson method converges beyond the contract range $b$. Increasing $b>0$ with continuation on $b$, our numerical tests show the steady convergence
for relatively large $b>0$.
    
\subsection{Numerical continuation on $b$}
We show the effect of numerical continuation \cite{K} on \eqref{MRE} through the composition of contractive maps by internal and external update. We start with a small $b > 0$. 
\begin{enumerate}
    \item External update: We update external variable $u$ with a given iterative scheme until the mapping is shown contractive. If fails, we terminate the algorithm.
    \item Internal update: We heuristically upgrade internal variable $b$ as a backward update and go back to the previous step.
\end{enumerate}
We compare the value of $b$ to see the effect of standard Banach iterate and the incorporation of Anderson type acceleration into Banach iterate.
We first let $n = 10$, and randomly generate a symmetric and positive definite (SPD) matrix $Q$. Then, we pick an initial $b>0$ such that contraction can happen. 
We now describe sensitivity for the coefficient $b$. For each given $b\in \mathbf{R}_{+}$ with different SPD $Q$ (even of the same dimension), the sensitivity range of $b$ for convergence will be different. So, $b$ is an important parameter to tune beyond the contraction constant. We start with a $b$ that is positive yet close to 0, so that the fixed point iterate will converge, to make sure that $b$ is able to grow up. We move $b$ forward so that external update works. Otherwise, we stop the entire algorithm. 

We also generalize from the case of $n = 10$ to $n = 200$ to see that Anderson acceleration methods work for the large scale Riccati equation in general.

We first discuss how the algorithm works when $u$ is an arbitrary squared size. We do initialization on a symmetric positive definite matrix $Q$ first to generate the starting iterate $u$. 
\begin{enumerate}
    \item We generate $Q$ as follows:
\begin{enumerate}
    \item We make $Q$ a symmetric matrix by $Q = \frac{1}{2} (Q + Q^{\top}).$
    \item  We compute the smallest eigenvalue of $Q$ and add appropriate regularization constant times identity to make sure $Q$ is a positive definite matrix. We make the minimum eigenvalue between 0 and 1. 
\end{enumerate}
    \item Initialization: With $u = 0$, fix $b \in \mathbf{R}_+$ such that the algorithm may converge, say $b = 0.05$. Now, we look at range of $b$ for which Anderson works. 
\end{enumerate}
    \noindent We now look at contraction mappings to update the external variable $u$. We first introduce the Banach iterate. 
    Fix $\epsilon>0$ as tolerance level, upgrade $u$ until $|u_{k+1} - u_k| < \epsilon.$ Let  
    \begin{equation} \label{Bana}
          u_{k+1} = -b\ u_k\ u_k + Q. 
    \end{equation}
    One would perform \eqref{Bana} with a batch of size $m$ to accelerate the convergence. 
    
\noindent Now, we observe the use of Nested Anderson map, with step size $m \ge 2$, from Section \ref{Nested}. We generate $u_1$ and $u_2$ with \eqref{Bana}. For $i = 1,2$, we store residual vectors $r_i = -b \ u_i \ u_i + Q - u_i$ and perform vectorization to transform $r_i$ from $n \times n$ matrix to a $n^2 \times 1$ vector, and we store $m$ of vectroized $r_i$ in matrix $r$ as column vectors. With $B_{ij} = (r_i, r_j)$, we perform \eqref{DIISSys} to determine weight and upgrade iterate by \eqref{ROMUpdate} as new basis elements. By taking a linear combination of past iterates, we update the current iterate. We implement successive fixed point iterates with Anderson depth of $2,3,4$ and printed out $|u_{k+1} - u_k|$ at each iterate to monitor acceleration convergence. At last, we solve for 4 mixing weights and form the final accelerated iterate as restart.
    

We list number of iterates with $|u_{k+1} - u_k|$ for different algorithms taking to converge on certain values of $b$ by making Table \ref{algo1} for standard  iterate with 3 steps, Table \ref{algo2.4} for Nested Anderson with 4 step, and Table \ref{algo2.7} for Nested 7-step Anderson. In each table, we first fix a symmetric positive definite matrix $Q$ of size $10 \times 10$, then we upgrade $b\in\mathbf{R}$ and we report the first $|u_{k+1} -u_k|$ below the tolerance level on the indicated iteration amount after the entire loop is finished. Anderson acceleration method works not only for small matrix of size $10 \times 10$ but also for the large scale Riccati equation in general. We illustrate this by considering $u$ as a matrix of $200 \times 200$ unknowns to solve, and acceleration methods works quite well, which shows the applicability of Anderson for large scale problems for general pattern in practices. Table \ref{ConvStep} shows that the Anderson method still works with the method of continuation for a large scale nonlinear system of equation when $u$ is of size $200 \times 200$. 

We stop as long as $|u_{k+1} -u_k|$ fails to converge to 0. We observe that the original fixed point iterate with large $b$  (say $0.5$) is never contractive. Employing Nested Anderson map makes contraction lasts longer, and Anderson does not require fixed point map contraction for the original function to hold. Therefore, we use the method of continuation that Anderson map accelerates contraction. 
    
    \begin{table}[h]
        \centering
        \begin{tabular}{|c|c|c|}
        \hline
            $b$ & Iterate & $|u_{k+1}-u_k|$\\
        \hline
            $0.05$ & 2 & $1.4918\times 10^{-10}$\\ 
        \hline
            $0.1$ & 11 & $8.7325\times 10^{-7}$\\ 
        \hline
        \end{tabular}
        \caption{Contraction continuation with standard Banach Iterate with 3 steps ($n=10$).}
        \label{algo1}
    \end{table}

    \begin{table}[h]
        \centering
        \begin{tabular}{|c|c|c|}
        \hline
            $b$ & Iterate & $|u_{k+1}-u_k|$\\
        \hline
            0.05 & 3 & $4.6539\times 10^{-7}$ \\ 
        \hline
            0.1 & 3 & $1.1879\times 10^{-7}$ \\ 
        \hline
        0.5 & 6 & $2.0847\times 10^{-7}$\\
        \hline
        1 & 6 & $6.9332\times 10^{-7}$ \\
        \hline
        1.5 & 8 & $8.9559\times 10^{-7}$ \\ 
        \hline
        \end{tabular}
        \caption{Contraction continuation with Nested 4 step Anderson acceleration ($n=10$).}
        \label{algo2.4}
    \end{table}

\begin{table}[h]
    \centering
    \begin{tabular}{|c|c|c|}
    \hline
        $b$ & Iterate & $|u_{k+1}-u_k|$ \\
    \hline
    1 & 4 & $1.4598\times 10^{-8}$ \\ 
    \hline
    2 & 5 & $1.0643\times 10^{-7}$ \\ 
    \hline
    5 & 7 & $8.2790\times 10^{-10}$ \\ 
    \hline
    10 & 9 & $1.7867\times 10^{-9}$\\ 
    \hline
    20 & 12 & $2.8622\times 10^{-8}$ \\ 
    \hline
    40 & 19 & $3.5862\times 10^{-9} $\\ 
    \hline
    45 & 18 & $7.8087\times 10^{-8}$\\ 
    \hline
    47.5 & 17 & $2.7625\times 10^{-9}$\\ 
    \hline
    \end{tabular}
    \caption{Contraction continuation with Nested 7 step Anderson acceleration ($n=10$).}
    \label{algo2.7}
\end{table} 

\begin{table}[h]
    \centering
    \begin{tabular}{|c|c|c|}
    \hline
    $b$ & Iterate & $|u_{k+1} -u_k|$\\
    \hline
        0.05 & 5 & $1.4765\times 10^{-8}$ \\ 
    \hline
        0.075 & 3 & $4.2009\times 10^{-7}$ \\ 
    \hline
        0.1 & 4  & $2.6738\times 10^{-7}$ \\ 
    \hline
        0.125 &  4 & $2.2680\times 10^{-8}$ \\ 
    \hline
        0.15 & 4 & $3.0553\times 10^{-7}$ \\
    \hline
        0.20 & 4 & $9.1989\times 10^{-7}$ \\
    \hline
       0.25 & 4 & $4.0601\times 10^{-7}$ \\
    \hline
       0.30 & 9 & $4.2750\times 10^{-7}$ \\ 
    \hline
    \end{tabular}
    \caption{Contraction continuation with Nested 7 step Anderson acceleration ($n=200$).}
    \label{ConvStep}
\end{table}

\begin{remark}
    Each time with a different $Q$, all tables would differ. A general pattern is that the larger $b$ becomes, the more iterate it will take to get to the desired level of convergence. We do regularization to guarantee the existence of weight coefficients, shown in Section \ref{extp}, if necessary.
\end{remark}


\begin{remark}
Since it is possible that the algorithm will converge to the desired tolerance level before the last step of iterate, it is possible to have redundant steps. We ignore the effects of redundant steps on the last iteration, whenever difference $|u_{k+1} -u_k|$ hits below the reaches the stopping criteria in some stage.      
\end{remark}

We bring summary of the section with the following observations with unknown matrix $u$ of size $10 \times 10$:
\begin{enumerate}
    \item Observation 1:  We use standard fixed point iterate: $\dot{u} = -b\ u\ u + Q.$ Small $b$ implies a fixed point contraction will happen, and Anderson acceleration help accelerate convergence fairly very much. For example, in the case when $b=0.1$, the standard case needs 9 more iterate than the Nested case.
    \item Observation 2: Test of Anderson for $b$ increment suggests that contraction can happen beyond the contraction range for the original equation.
    \item Observation 3: Given a large $b$, Anderson stop converging, therefore we adapt the idea of continuation: we enlarge $b$ sequentially, then we apply acceleration. As shown in Table \ref{algo2.4} and Table \ref{algo2.7}, Anderson shows contraction with continuation further. 
\end{enumerate}
At last, we refined matrix $u$ of size $200 \times 200$ unknowns for Riccati model and the same observation from holds as shown in Table \ref{ConvStep}, which illustrates that Anderson map works for a nonlinear system in general. Therefore, we have remedied contraction map with Anderson to make contraction work for nonlinear system. 

In summary, Anderson and method of continuation work for nonlinear systems. We may use a different method to remedy small contraction range of $b$ when using the pre-conditioned fixed point method. 

\section{Saddle point problems}\label{Exper4}
In this section, we focus on large scale semi ill-posed linear system of equations $Ax = b$, specifically the Saddle Point Problem, with a given $n \in \mathbf{Z}_+$, and $A \in \mathbf{R}^{n \times n}$, and $b \in \mathbf{R}^n$. We solve it by Cauchy's method and its variant, along with ROM, and Sampled ROM. We report the relative residual for each algorithm. We specify Stokes-like systems as an example for Saddle Point Problem. 
\subsection{General description}\label{Saddle}
In this section we discuss the saddle point problems.
Consider the constrained minimization
$$
\min\quad F(x) \mbox{  subject to   }E(x )= b
$$
Define the Lagrangian
$$
L(x, p) = F (x) + ( E(x)-b, p).
$$
It follows from the Lagrange multiplier that the pair $(x, p)$ satisfies the saddle-point problem
$$
{\cal A}(x,p) =\left(\begin{array}{cc} F^\prime (x) +
(E^\prime(x),p)\\ E(x)-b +
 \epsilon p \end{array}\right)= \left(\begin{array}{cc}
    0\\
    0 
 \end{array}\right),
$$
where we assume $dE: T_xX \to T_{E(x)}Y$ is surjective and $\epsilon > 0$. 
It is a skew form and
${\cal A}$ is nonnegative 
definite if $(F^\prime(x)\phi,\phi)>0$ for all $\phi \in X.$

\noindent \underline{\textbf{Example}: (Linear saddle point problem)} Given $D$ a symmetric positive definite matrix, if $F(x)=(x,Dx)-(a,x)$ and $E(x) = Ex-b$,
then we have 
$$
{\cal A}(x,p) =\left(\begin{array}{cc} D& E^* \\ \\  -E & O\end{array}\right) \left(\begin{array}{cc} x \\ \\p\end{array}\right) = A \left(\begin{array}{cc} x \\ \\p\end{array}\right) = \left(\begin{array}{cc} a \\ \\b\end{array}\right).
$$

\subsection{Saddle point problem 1}\label{S1} 
Recall from Section \ref{Saddle} in the case of  two dimensinal Stokes system for $x=(u,v)$  given $D$ as a symmetric positive definite matrix, we consider the saddle point system given by a specified $(D,E)$ as a constrained minimization problem:
\begin{equation} 
 {\frac{1}{2}(u,Du) + \half (v,Dv) -(\vec{f},u)} \mbox{ subject to  }
 Ex=E_1u+E_2v=g.
\end{equation}
We consider a divergence free case, i.e., the case when $g = 0$. In addition, for simplification, we treat $v=0$. The necessary optimality condition implies $A(u,v,p)= b$:
\begin{equation}\label{Sad1}
\underbrace{\left[ 
\begin{array}{c|c} 
  D & E^{\top} \\ 
  \hline
  E & \mathbf{0} 
\end{array} 
\right]}_{A} \underbrace{\begin{bmatrix}
    u\\
    p
\end{bmatrix}}_{x} = \underbrace{\begin{bmatrix}
    f\\
    0
\end{bmatrix}}_{b}.\end{equation} 
We consider homogeneous boundary condition for $u$. Let $n \in \mathbf{Z}_+$ be the dimension of for 2nd order central difference matrix $T_n$ given by
\begin{equation}\label{2nddiff}
T_n = \begin{bmatrix}
2 & -1 & 0 & \cdots & 0 \\
-1 & 2 & -1 & \cdots & 0 \\
0 & -1 & 2 & \cdots & 0 \\
\vdots & \vdots & \vdots & \ddots & -1 \\
0 & 0 & 0 & -1 & 2
\end{bmatrix}.
\end{equation}
We use the homogeneous boundary condition on $(u,v)$to define $(D,E)$. We first model $D$ by:
\[
\mathbf{D} = I_n\otimes T_n + T_n \otimes I_n,
\]
where $I_n$ is the identity matrix, and $T_n$, defined by \eqref{2nddiff}, is the central difference approximation of $-u_{xx}$ with homogeneous boundary condition on $u$ of size $n^2 \times n^2$. $D$ is a positive definite matrix with $\min(\text{eig}(D)) = 0.1620$. Here $D$ gives energy for minimization, and $E$ gives the divergence free constraint 
$Ex = (u_x+v_y)$ for $u_x$ with homogeneous boundary condition on $u$. 

Now, we construct $E$ by: 
\[
\mathbf{E} = I_n \otimes C_n + C_n \otimes I_n,
\]
where we set the first order central difference matrix $C_n$ of size $n^2 \times n^2$ by:
\begin{equation} \label{1stdiff}
C_n 
= \begin{bmatrix}
0 & 1 & 0 & \cdots & 0 \\
-1 & 0 & 1 & \cdots & 0 \\
0 & -1 & 0 & \ddots & 0 \\
\vdots & \ddots & \ddots & \ddots & 1 \\
0 & \cdots & 0 & -1 & 0
\end{bmatrix}.
\end{equation}
That is, given different condition with appropriately defined the corresponding $(D,E)$.
For $n = 10$, we have the matrix $A$ of size $2n^2 \times 2n^2$ in \eqref{Sad1}. Since $\min(\text{eig}(A)) <0$ and $\max(\text{eig}(A)) > 0$, $A$ is an indefinite matrix.

We discuss experiment setup. We set the regularization constant $10^{-30}$ and $m = 15$ as the amount of loop to go over. For ROM, we store $m$ residual vectors into a residual matrix, and perform Anderson acceleration; For sampled ROM, we store $m$ vectors into a residual matrix to perform the first round of Anderson acceleration, and their resulting performance $m$ vectors into another residual vector to perform another round of Anderson acceleration.

\begin{figure}
    \centering
    \includegraphics[width=0.75\linewidth]{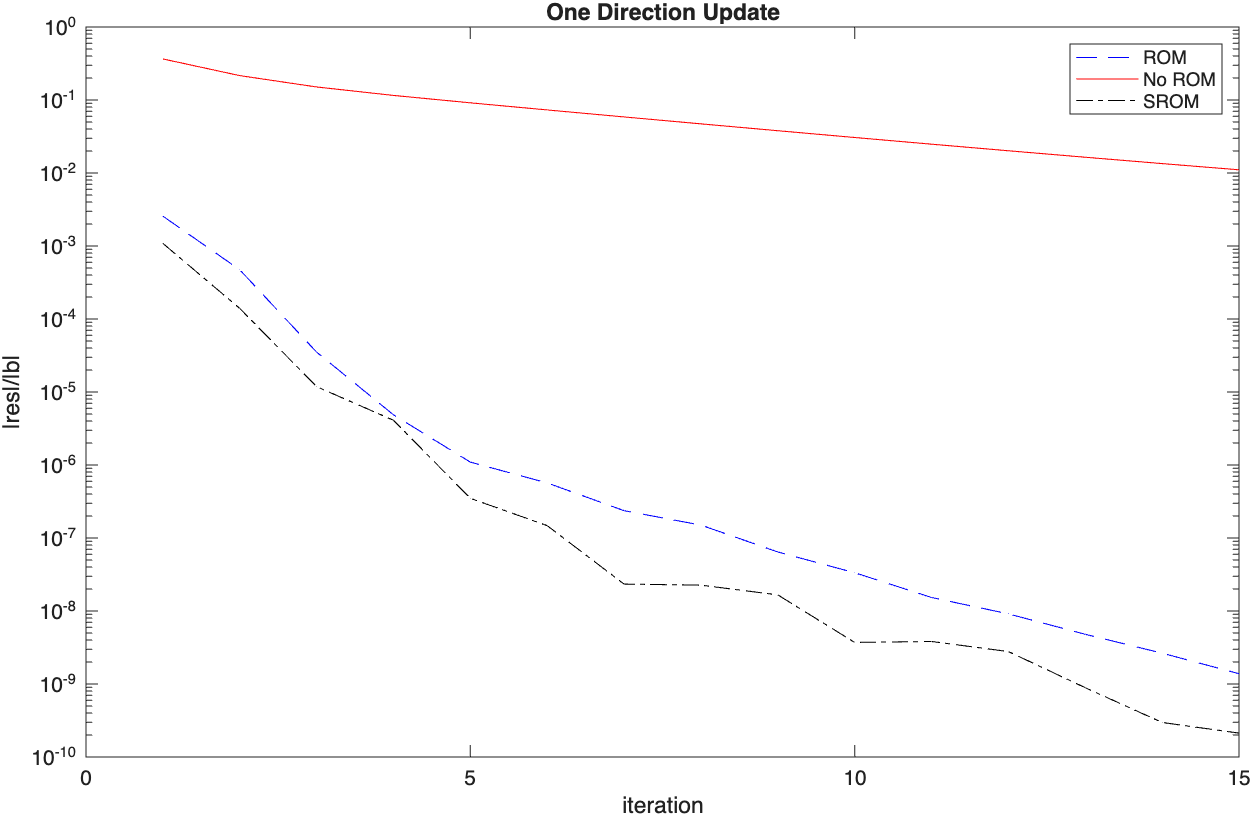}
    \caption{Saddle Point Problem 1 (One direction update: $r = b-Ax$): comparison of standard case, accelerated case, and sampled case. }
    \label{StokesUno} 
\end{figure}

\begin{figure}
    \centering
    \includegraphics[width=0.75\linewidth]{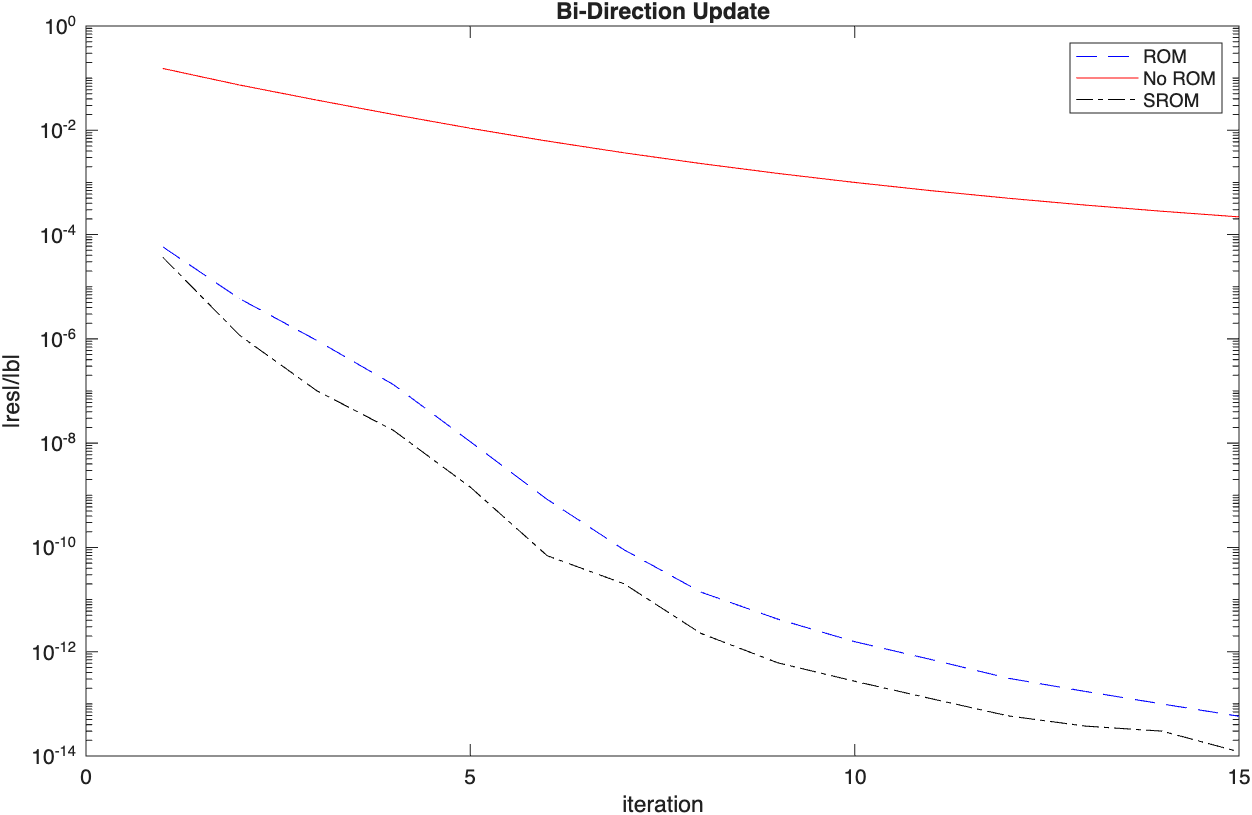}
    \caption{Saddle Point Problem 1 (Two directions' update: $r = b-Ax$, $Ar$): comparison of standard case, accelerated case, and sampled case.}
    \label{StokesDos} 
\end{figure}

We plot relative residual for all iterate and we shall see that: ROM accelerates convergence, when performing a one-direction update (c.f. Figure \ref{StokesUno}), or a two direction update (c.f. Figure \ref{StokesDos}). As suggested from pictures, we see that further improvement of ROM can be done through the Sampled ROM. In addition, the performance with a two directions' update as a multi step subspace update will significantly improve the performance of a one direction update. 


\begin{remark}
For general boundary condition (e.g., von-Nuemann condition), we adjust $D$ and $E$ appropriately; With $D$ and $C_n$ same as Example above and with $E_1 = C_n \otimes I_n, E_2 = I_n \otimes C_n$, our experiment found that same conclusion from the example still holds.    
\end{remark}
 

\subsection{Saddle point problem 2}\label{Experi5}
In this section, we use a variant of Cauchy step size update for solving $Ax = b$. We generate an approximated matrix $\bar{A}$ for $A$. In one direction update, we use an alternative step size \eqref{ACUno} than Cauchy step \eqref{CUno} as a variable step size update. By using $\bar{A}$, we sacrifice optimality in each step in exchange for efficiency. Cauchy method \eqref{CUno} gives optimal step size of the steepest descent, and \eqref{ACUno} gives a cheaper surrogate that approximates optimal step size of descent direction but allows for faster computation with cost reduction. We iterate with $\bar{A}$ and residual from the original equation $r = b-Ax$ to compute variable step size update in \eqref{ACUno}, then correct with a ROM step that uses a richer space of prior residuals to better approximate the true solution trajectory. For example, in saddle-point systems (like Stokes or Navier–Stokes), solving with $A$ is hard due to its null space or ill-conditioning. 

We shall treat $\bar{A}$ as a preconditioner of $A$, but we make use of $\bar{A}$ as a preconditioner in an unconventional way without inverting matrix $\bar{A}$. We use $\bar{A}$ to approximate the action of $A$ in computing variable step size update, which improves computational efficiency in place of exact information from $A$. In other words, we have applied a cheap surrogate for $A$ to implicitly modify the step size of a given search direction, which matches the goal of applying a preconditioner with cheaper replacement for $A$ to modify the search direction implicitly. We use Anderson to accelerate our proposed method. 

In terms of two directions' subspace update implementation, one considers \eqref{ACDos} than the Cauchy step \eqref{CDos}. With exactly the same argument from one direction update, we find that the same conclusion from 1 direction subspace update also hold for 2 directions subspace update.



Consider the saddle point problem
\begin{equation}
A = \left(\begin{array}{cc} D+\delta F & E^* \\ \\ E & \huge{0} \end{array} \right),    
\end{equation}
with the same $D$, $E$ as in Section \ref{S1}. We set the perturbation $\delta F$ on $D$ by: 
For example,
$$
\delta F = \mathbf{1}^\top E, $$
where $\mathbf{1}$ is the matrix of all ones.

For $\bar{A}$, we have: 
\begin{equation}
    \bar{A} =  \left(\begin{array}{cc} D & E^* \\ \\ E & 0 \end{array} \right).
\end{equation}
Here, $A$ is a low rank perturbation of $\bar{A}$, it is easier to work with $\bar{A}$ to reduce computational cost. We shall see that the update with \eqref{ACUno} and \eqref{ACDos} are nearly optimal due to low rank perturbation between $A$ and $\bar{A}$. We set the regularization constant as $10^{-22}$ and $m = 15$.

\begin{figure}
    \centering
    \includegraphics[width=0.75\linewidth]{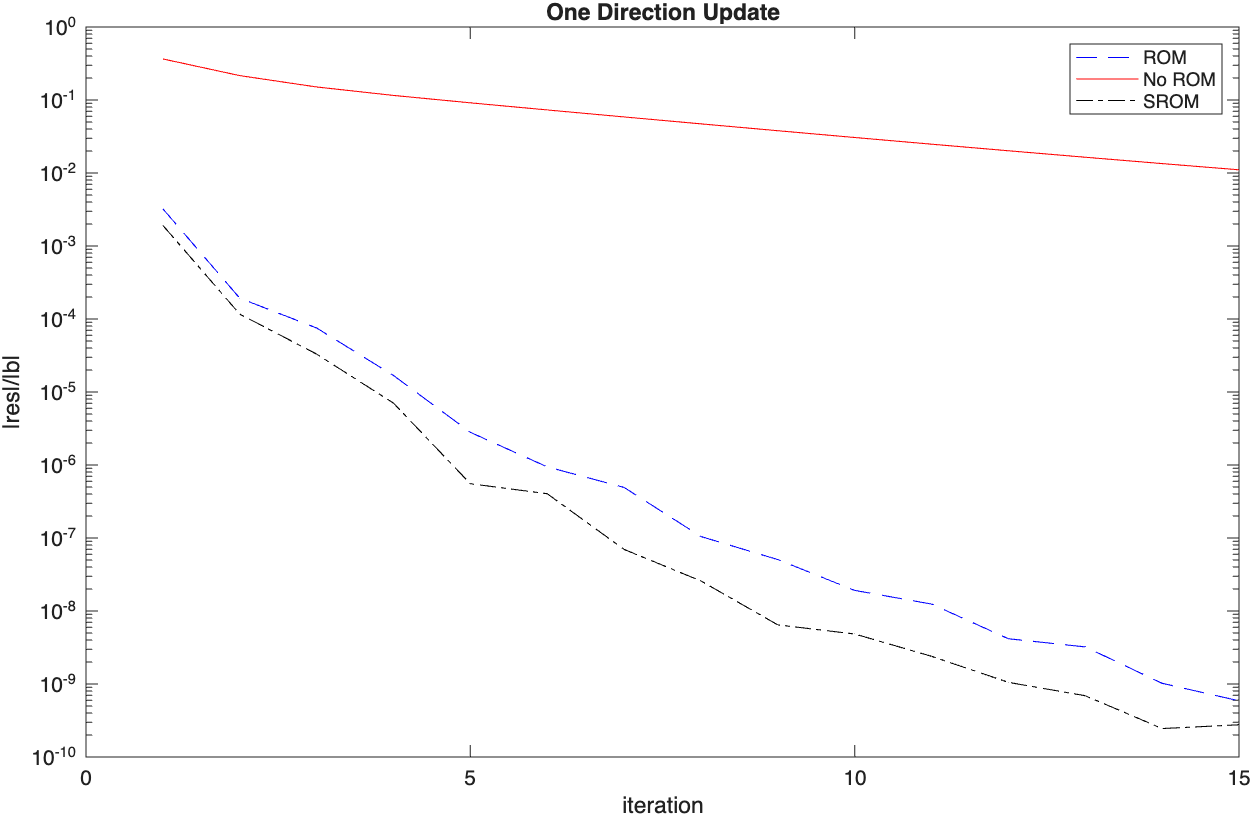} 
    \caption{Saddle Point Problem 2 (One direction update $r=b-Ax$): comparison of standard case, accelerated case, and sampled case.} 
    \label{PrecondtioningEx2Uno}
\end{figure}

\begin{figure}
    \centering
    \includegraphics[width=0.75\linewidth]{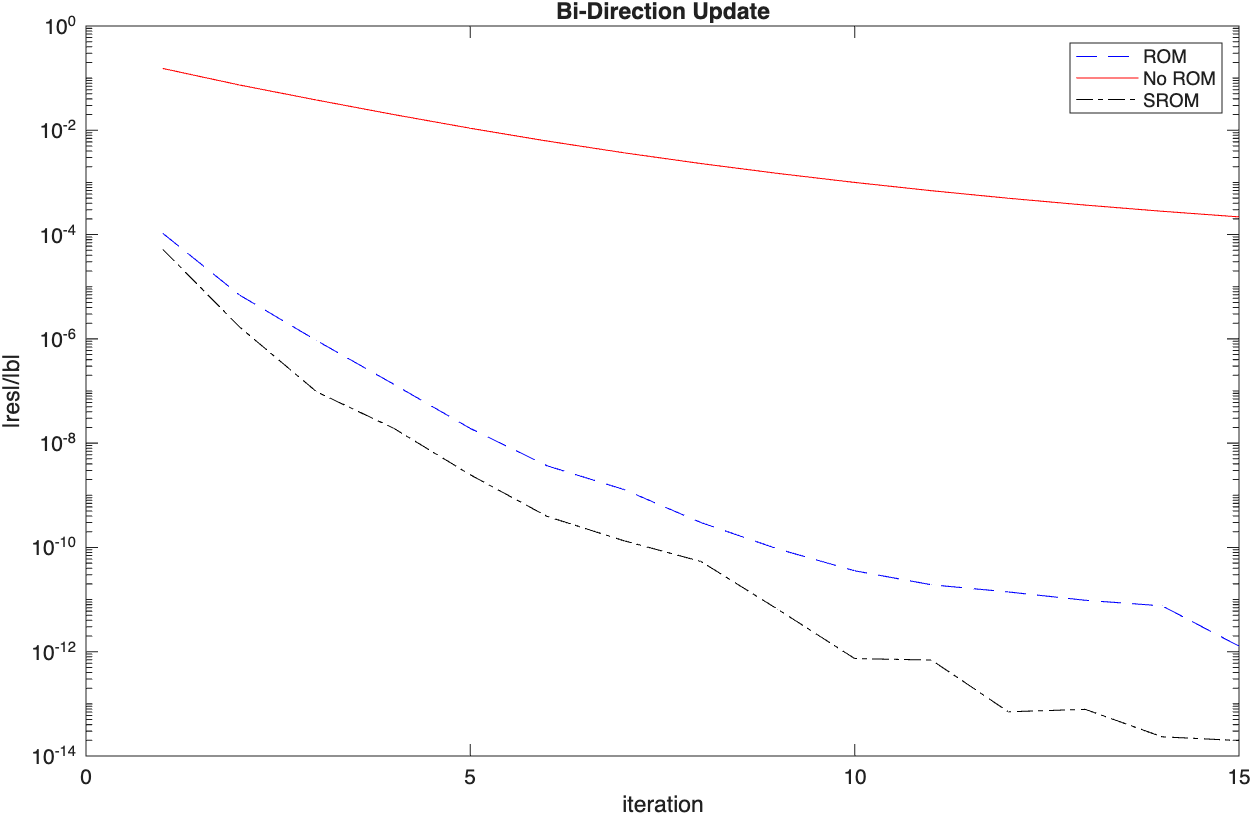}
    \caption{Saddle Point Problem 2 (Two directions' update $r=b-Ax$, $\bar{A}r$): comparison of standard case, accelerated case, and sampled case.} 
    \label{PrecondtioningEx2Dos}
\end{figure}

We plot relative residual and we shall see that: ROM speeds up the convergence, when performing a one-direction update (c.f. Figure \ref{PrecondtioningEx2Uno}), or a two direction update (c.f. Figure \ref{PrecondtioningEx2Dos}). Further improvement can be done through the Sampled ROM. In addition,  the performance with a two directions' update will boost the performance of a one direction update. 

\begin{remark}
    Similar to example in Section \ref{S1}, when constrained with a boundary condition (e.g., Dirichlet condition), we adjust $D$ and $E$ appropriately.
\end{remark}

\section{Conclusions}
\label{sec:conclusions}
In this work, we have developed and investigated iterative schemes based on Reduced Order Method to generate basis and accelerate algorithms by refining iterative scheme by developing variable step size. At the end, we have tested algorithms in applications for both linear and nonlinear systems.



\section*{Acknowledgments}
T. Xue would like to thank Dr. Bangti Jin and Dr. Fuqun Han in giving him suggestions to revise the paper.


\end{document}